\documentclass[12pt,letterpaper,spanish]{article}
\usepackage[utf8x]{inputenc}
\usepackage{ucs}

\usepackage{amsfonts,amssymb,amstext,amsmath,amscd,latexsym,amsthm}

\usepackage{pgf,tikz} 
\usetikzlibrary{arrows}

\usepackage[spanish]{babel}
\usepackage{hyperref}

\usepackage{multicol}
\usepackage{enumerate}
\usepackage{answers}
\Newassociation{sol}{Solution}{ans}

\usepackage{slashbox}
\usepackage{float}       
\usepackage{color}

\theoremstyle{plain}
\theoremstyle{definition}

\newtheorem{remark}{Observaci\'on}

\title{Una relación entre la\\
Geometría y el Algebra\\
(Programa de Erlangen)
}
\author{José Ricardo Arteaga B.\\
Universidad de los Andes,\\
Departamento de Matem\'aticas\\
\texttt{jarteaga@uniandes.edu.co}}
\date{Junio 2011}

\usepackage{fancyhdr,lastpage}
\begin{document}
\maketitle

\begin{abstract}
The three key documents for  study  geometry are: 1) `` The Elements'' of Euclid, 2) the  lecture by B. Riemann at Göttingen in 1854 entitled ``Über die Hypothesen welche der Geometrie zu Grunde liegen" (On the hypotheses which underlie geometry) and 3) the ``'Erlangen Program", a document written by F. Klein (1872) on his income as  professor at the Faculty of Philosophy and the Senate of the Erlangen University.
The latter document F. Klein introduces the concept of group as a tool to study geometry. The concept of a group of transformations of space  was known at the time. \ \
The purpose of this informative paper is to show a relationship  between geometry and algebra through an example,  the projective plane. Erlangen program until today continues being a guideline of how to study geometry.
\end{abstract}

\begin{abstract}
Los tres documentos fundamentales para el estudio de la Geometría son: 1) ``Los Elementos'' de Euclides, 2) la conferencia de B. Riemann ``Sobre las hipótesis que están los fundamentos de la Geometría'' (1854) con motivo de su habilitación para ser profesor universitario (Universidad de G\"ottingen) y 3) el ``Programa de Erlangen'', documento escrito por  F. Klein  (1872) con motivo de su ingreso  como profesor a la Facultad de Filosofía y al Senado de la Universidad de Erlangen.
En este último documento F.Klein introduce el concepto de grupo como una herramienta para estudiar  Geometría. El concepto de grupo de transformaciones de un espacio ya era conocido en ese entonces. \\
El objetivo de este documento divulgativo es mostrar una relación de la Geometría y el Algebra tomando como ejemplo   el plano proyectivo.  El programa de Erlangen sigue marcando hasta hoy día una directriz de cómo estudiar y hacer Geometría moderna. 
\end{abstract}

\section*{Introducción}
La historia de los orígenes de la Geometría es igual a los de la Aritmética. Los conceptos geométricos más antiguos provienen de diversas culturas y son una consecuencia de las actividades prácticas del hombre. 
En el siglo VII AC la geometría pasó de Egipto a Grecia. La escuela pitagórica, secta religiosa-filosófica idealista, hizo grandes contribuciones a la geometría que se plasmaron luego en un documento del siglo III AC conocido como ``Los Elementos'' de Euclides.\footnote{Euclides de Alejandría (325-265 AC, Egipto) escribió ``Los Elementos". La primera versión impresa  traducida del árabe al latín  de Los Elementos apareció en 1482.}   En este trabajo, la geometría se presentó como un sistema muy sólido, y sus fundamentos no sufrieron cambios esenciales hasta la llegada de N.I. Lobatchevski.\footnote{Nikolái Ivánovich Lobatchevski (1792-1856, Rusia), profesor de la Universidad de Kazán, Rusia, el 23 de febrero de 1826 leyó una memoria sobre la teoría de las paralelas en una sesión de la facultad de Física y Matemáticas, y en 1829 publicó su contenido en la revista de la Universidad de Kazán.\cite{Alexandrov} } Lobatchevski escribió en 1835 ``Los nuevos Elementos de Geometría'', con el cual se dió inicio a una nueva geometría no euclidiana conocida geometría de Lobatchevski o geometría hiperbólica.

El 10 de junio de 1854, B. Riemann\footnote{Georg Friedrich Bernhard Riemann (1826-1866), Alemania} dió una conferencia en la Universidad de Gotinga para ser profesor universitario. El tema de la conferencia por sugerencia de Gauss, su protector y antiguo profesor durante la licenciatura y el doctorado, fue sobre Geometría, y se tituló ``Sobre las hipótesis que están en los fundamentos de la geometría''. 

En 1872 F. Klein \footnote{Felix Christian Klein (1849-1925), Alemania} presentó un programa de investigación con motivo de su ingreso  como profesor a la Facultad de Filosofía y al Senado de la Universidad de Erlangen, Alemania. En este trabajo Klein propuso una nueva solución al problema de cómo clasificar y caracterizar las geometrías existentes sobre la base de la geometría proyectiva y la teoría de grupos.

\begin{figure}[H]
\begin{minipage}[t]{5cm}
\includegraphics[scale=0.5]{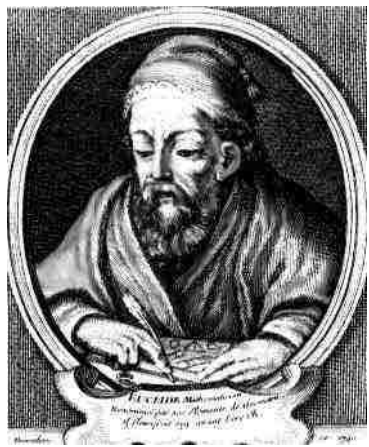}
\end{minipage}
\hfill\begin{minipage}[t]{5cm}
\includegraphics[scale=0.5]{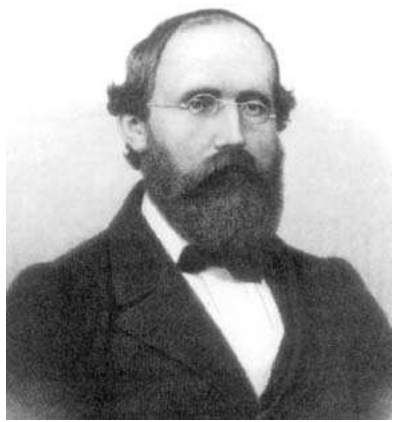}
\end{minipage}
\hfill\begin{minipage}[t]{5cm}
\includegraphics[scale=0.4]{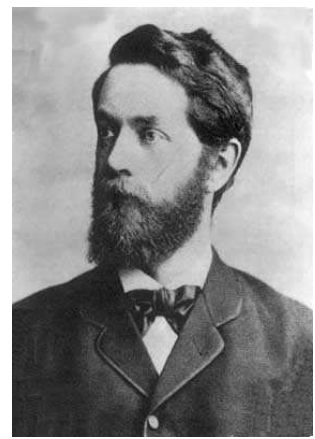}
\end{minipage}
\caption{Euclides -- B. Riemann -- F. Klein}
\end{figure}

\subsection*{La perspectiva}
El problema inicial es, ¿cómo dibujar en el papel algo que ven mis ojos? El primer paso en tratar de encontrar leyes que resolvieran este problema lo dieron los pintores. Por ejemplo el artista Leon Battista Alberti (1404-1472) y  el ingeniero y arquitecto Florentino Fillipo Brunelleschi (1377-1446)  desarrollaron métodos y una teoría matemática de la perspectiva la cual fue seguida por los pintores de renacimiento \cite{eb}.

\begin{figure}[H]
\includegraphics[scale=0.5]{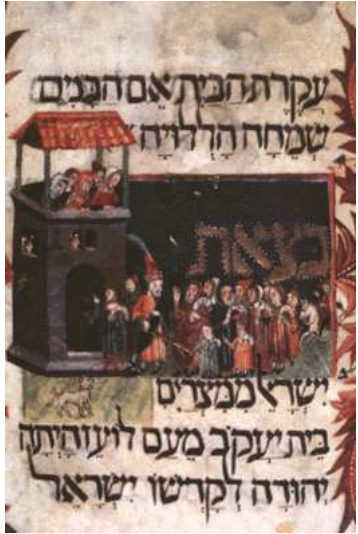}
\caption{Pintura antes de la perspectiva. 
Panel inicial del Salmo de la Hagadá de Kaufmann. España finales del siglo  XIV. Observe que las lineas del balcón parecen no converger.}
\end{figure}

\begin{figure}[H]
\begin{center}
\includegraphics[scale=0.4]{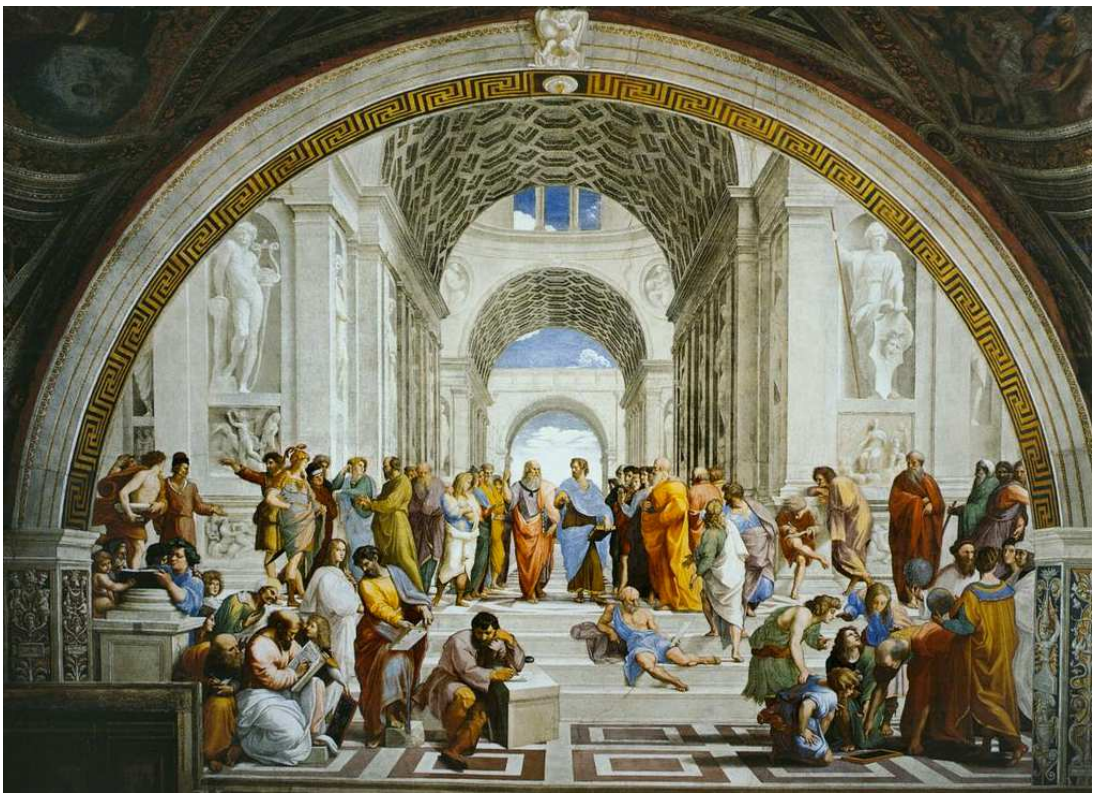} 
\caption{Pintura después de la perspectiva.
Rafael Sanzio (1483-1520), Escuela de Atenas. Observe que las líneas superiores de los muros verticales parecen converger a un punto.}
\label{Fig:fig2-con_perspectiva}
\end{center}
\end{figure}

Para comprender mejor un dibujo en perspectiva primero veamos la proyección de un objeto sobre un plano. En el espacio tridimensional $\mathbb{R}^{3}$ consideremos un punto $O$ fijo (sobre un plano $\beta$) que llamaremos centro de la proyección y un plano fijo $\alpha$, paralelo a $\beta$, llamado plano de la proyección, que no contiene a $O$. La proyección sobre $\alpha$ de un punto $A$ en el espacio $\mathbb{R}^{3}$, respecto a $O$, es el punto $\tilde{A}$ tal que $A$, $O$ y $\tilde{A}$ son colineales.

\begin{figure}
\begin{center}
\includegraphics[scale=0.5]{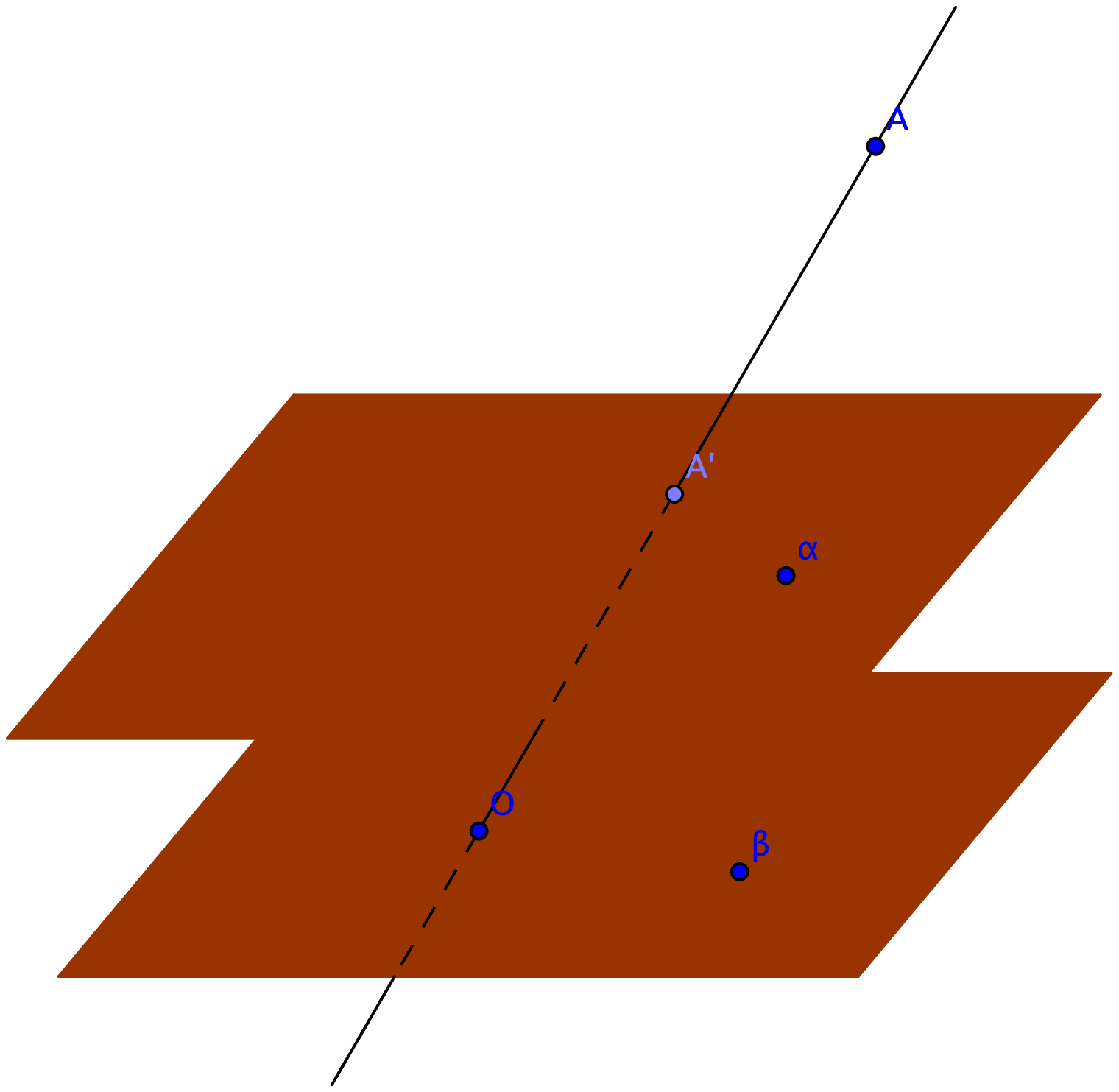}
\caption{Proyección de $A$ en $\tilde{A} \in \alpha$ con centro de la proyección en $O \in \beta$. Dibujo  hecho con GEOGEBRA \href{http://www.geogebra.org/cms/}{http://www.geogebra.org/cms/}}
\end{center}
\end{figure}

\begin{remark}
Debe ser claro que,
\begin{enumerate}
\item
$\tilde{A}$ es la imagen (proyección) en $\alpha$ no de un único punto, si no de todos los puntos sobre la recta $OA$. Es decir dos puntos $A_{1}$, $A_{2}$ diferentes a $O$ en $\mathbb{R}^{3}$ tienen la misma imagen (proyección) en $\alpha$ si y solamente si $A_{1}$, $A_{2}$ y $O$ son colineales.
\item
existen puntos en $\mathbb{R}^{3}$ que no tienen imagen en $\alpha$. Todos aquellos puntos que pertenecen al plano $\beta$ paralelo a $\alpha$ y que contiene a $O$.
\item
si tomamos en el sistema coordenado cartesiano $(x,y,z)$ en $\mathbb{R}^{3}$ el plano $\alpha$ como el plano $z=1$, el centro de la proyección $O$ como el origen del sistema, entonces las coordenadas de un punto $\tilde{A}$ son $(\tilde{x}, \tilde{y}, 1)$. En este caso si $A(x,y,z)$ no pertenece al plano $\beta$ ($z=0$), se proyecta en $\tilde{A}(\tilde{x}, \tilde{y}, 1)$ donde
\begin{equation}
\tilde{x} = \dfrac{x}{z}, \quad \tilde{y} = \dfrac{y}{z}
\end{equation}
\end{enumerate}
\end{remark}

\subsection*{Qué es Geometría?}

Podríamos decir que Geometría es el estudio de objetos que están en cierto ambiente (espacio). El estudio se refiere a la búsqueda de invariantes de relaciones  entre  los objetos (independientes del sistema de referencia) bajo la acción de un grupo. Para estudiar un tipo de Geometría debemos entonces definir:
\begin{enumerate}
\item 
Un espacio ambiente.
\item 
Los objetos que vamos a estudiar en este espacio.
\item 
Un grupo que define los movimientos. Los movimientos son transformaciones del espacio que preservan ciertas relaciones entre los objetos, las cuales estamos interesados en estudiar. Estas transformaciones forman un grupo.
\end{enumerate}
 
Por ejemplo, la geometría euclidiana en el plano estudiada en la escuela y el bachillerato, el espacio ambiente es el plano, los objetos son puntos, rectas, polígonos, circunferencias, etc. y el grupo de transformaciones es el grupo de transformaciones del plano que preservan la distancia (grupo euclidiano).

\section*{El plano proyectivo $\mathbb{RP}^{2}$}

El plano proyectivo es una extensión del concepto de plano euclidiano que conocemos.  En el plano proyectivo cualquier par de rectas se intersectan a diferencia de las rectas en el plano euclidiano. Veamos el siguiente modelo del plano proyectivo.

\subsection*{Modelo del plano proyectivo $ \mathbb{RP}^{2}= \mathbb{R}^{3}_{0}/ \sim $}

En el espacio tridimensional sin el origen $ \mathbf{R}^{3}_{0} = \mathbf{R}^{3} \backslash \{(0,0,0)\} $  consideremos  la siguiente relación de equivalencia: Dados dos puntos $A_{1} = (a,b,c)$ y $A_{2} = (d,e,f)$, $A_{1}$ es equivalente con $A_{2}$, $A_{1} \sim A_{2}$, si y solamente si existe un número real $k\neq 0$, tal que $a= kd$, $b=ke$, y $c=kf$,
\begin{equation}\label{Eq:relacion_equiv}
\quad (a,b,c) \sim (d,e,f) \iff (a,b,c) = k (d,e,f) , \text{para algún } k \in \mathbb{R} \backslash 0.
\end{equation}
El plano proyectivo lo definimos como $\mathbb{RP}^{2} = \mathbb{R}^{3}_{0}/ \sim $.
Los objetos fundamentales del plano proyectivo son puntos y rectas:
\begin{enumerate}
\item \textbf{Puntos.}
Los puntos proyectivos del plano proyectivo son las clases de equivalencia de la relación $\sim$. Si $A=[(x_{1}, x_{2}, x_{3})]$, donde los corchetes significan clase de equivalencia, entonces escribimos  $ [x_{1}: x_{2}: x_{3}] $ y se llaman coordenadas homogéneas del punto $A$ del plano proyectivo. En este modelo podemos ``visualizar'' un punto proyectivo como una recta en el espacio tridimensional $\mathbb{R}^{3}$ que pasa por el origen.
\item  \textbf{Rectas.} Dos puntos proyectivos diferentes en el plano proyectivo determinan una única recta proyectiva. En este modelo una recta proyectiva  del plano proyectivo es un plano en $ \mathbb{R}^{3}_{0} $ que pasa por el origen.
\end{enumerate}

Este proceso descrito anteriormente es la proyectivización de $\mathbb{R}^{3}$, es decir el plano proyectivo es la proyectivización de $\mathbb{R}^{3}$. Similarmente  la recta proyectiva es la proyectivización de $\mathbb{R}^{2}$.

\subsection*{El grupo proyectivo $PGL(3) = GL(3) / \sim$}

El grupo general lineal $ GL(3) $ es el grupo de matrices $ 3\times 3 $ no singulares,
\begin{equation}
GL(3) :=
\left\{
g = \left(
\begin{matrix}
a_{11} & a_{12} & a_{13} \\
a_{21} & a_{22} & a_{23} \\
a_{31} & a_{32} & a_{33} \\
\end{matrix}
\right)
\mid
\det (g) \neq 0
\right\}.
\end{equation}

Una matriz $ [a_{ij}] \in GL(3) $ se llama matriz escalar si $ a_{ii} = \lambda $ y $a_{ij}=0$ ($i\neq j$, $ \lambda \neq 0 $). El conjunto de matrices escalares con la multiplicación corriente de matrices forman un grupo. Denominemos este grupo por $ S $.
\begin{equation}
S :=
\left\{
s = \left(
\begin{matrix}
\lambda & 0 & 0 \\
0& \lambda & 0 \\
0 & 0 & \lambda\\
\end{matrix}
\right)
\mid
\lambda \neq 0
\right\}.
\end{equation}
$S$ es el centro de $GL(3)$, es un subgrupo normal y abeliano.

Definamos la siguiente relación de equivalencia: Dos matrices $g_{1}, g_{2} \in GL(3)$ son equivalentes, $g_{1} \sim g_{2}$ si existe una matriz escalar $s\in S$, tal que $g_{1} = s g_{2}$,
\begin{equation}
g_{1}, g_{2} \in GL(3), \quad g_{1} \sim g_{2} \iff g_{1} = s g_{2}, \text{ para alguna matriz } s\in S
\end{equation}

El grupo proyectivo $PGL(3)$ se define como $GL(3) / \sim $. Sus elementos son clases de equivalencia, es decir que un elemento $a$ de $PGL(3)$ es 
\begin{equation}
a=[g] = \{\lambda g \mid g\in GL(3), g \text{ fijo }, \lambda \in \mathbb{R} \backslash \{0\}\}
\end{equation}

\subsection*{Transformaciones}
Las transformaciones proyectivas (homografías) del plano proyectivo son el resultado de la acción del grupo proyectivo $PGL(3)$ sobre los puntos del plano proyectivo. Un movimiento de un punto $ A $ en el plano proyectivo lo podemos definir como $ [g] \cdot A$, donde $ [g] $ es un elemento del grupo proyectivo.
\begin{equation}
\begin{split}
A=&[x_{1}: x_{2}: x_{3}]\in \mathbb{RP}^{2} \\ 
[g] \in PGL(3), \quad [g] &=
\left[\left(
\begin{matrix}
a_{11} & a_{12} & a_{13} \\
a_{21} & a_{22} & a_{23} \\
a_{31} & a_{32} & a_{33} \\
\end{matrix}
\right)
\right]\\
\text{ la acción es: } [g]\cdot A &=
\left[\left(
\begin{matrix}
a_{11} & a_{12} & a_{13} \\
a_{21} & a_{22} & a_{23} \\
a_{31} & a_{32} & a_{33} \\
\end{matrix}
\right)
\right]
[x_{1}: x_{2}: x_{3}]^{T}
= [y_{1}: y_{2}: y_{3}]^{T}
\end{split}
\end{equation}

Podemos dividir los puntos del plano proyectivo en dos clases, los puntos propios y los impropios. Los propios son aquellos que su tercera coordenadas homogénea es diferente de cero y los impropios los que tienen su tercera coordenada homogénea igual a cero.  Note que una homografía puede llevar un punto propio a uno impropio y visceversa.

\subsection*{Invariantes}

En la recta proyectiva un punto propio tiene coordenadas homogéneas $[x_{1}: x_{2}] = [a: 1]$ donde $a = x_{1}/ x_{2}$ es un número real. En este caso denotaremos por $A(a)$ el punto con su coordenada.

La \textit{razón doble} entre cuatro puntos propios distintos $A=[a_{1}:a_{2}]$, $B=[b_{1}: b_{2}]$, $C=[c_{1}: c_{2}]$ y $D= [d_{1}: d_{2}]$ de la recta proyectiva $\mathbb{RP}^{1}$ se define como 
\begin{equation}
(A,B;C,D) = 
\dfrac{
\left|
\begin{matrix}
a_{1} & a_{2}\\
c_{1} & c_{2}\
\end{matrix}
\right|
\left|
\begin{matrix}
b_{1} & b_{2}\\
d_{1} & d_{2}\
\end{matrix}
\right|
}{
\left|
\begin{matrix}
a_{1} & a_{2}\\
d_{1} & d_{2}\
\end{matrix}
\right|
\left|
\begin{matrix}
b_{1} & b_{2}\\
c_{1} & c_{2}\
\end{matrix}
\right|
}
\end{equation}En el caso de puntos propios distintos sobre una recta la razón doble se simplifica a
\begin{equation}
(A,B;C,D) = \dfrac{c-a}{c-b}\cdot \dfrac{d-b}{d-a}
\end{equation}
El único invariante numérico de la recta proyectiva $\mathbb{RP}^{1}$ bajo la acción del grupo proyectivo $PGL(2)$ es la razón doble.

\begin{figure}
\begin{center}
\includegraphics[scale=0.5]{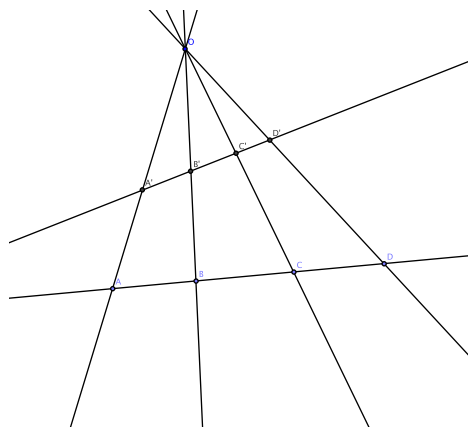}
\caption{Las razones dobles son iguales, $(A,B;C,D) = (A',B';C',D') = 1.35$. Dibujo y cálculos hechos con GEOGEBRA \href{http://www.geogebra.org/cms/}{http://www.geogebra.org/cms/}}
\end{center}
\end{figure}

\section*{Conclusiones}

\begin{itemize}
\item
Se llama geometría proyectiva a la pareja $(\mathbb{RP}^{n}, PGL(n))$, donde $\mathbb{RP}^{n}$ es la proyectivización de $\mathbb{R}^{n}$ y $PGL(n)$ es el grupo proyectivo.\footnote{$(\mathbb{RP}^{1}, PGL(1))$ es la recta proyectiva, $(\mathbb{RP}^{2}, PGL(2))$ es el plano proyectivo, $(\mathbb{RP}^{3}, PGL(3))$ espacio proyectivo, etc.}. En general,  los puntos proyectivos  son líneas rectas en $\mathbb{R}^{n}$ que pasan por el origen $O$. 
\item
Imaginemos  que estamos observando un objeto con un solo ojo desde un único punto $O$. Entonces todo lo que vemos es gracias a los rayos de luz que entran a nuestro ojo. Estas son líneas rectas o puntos proyectivos.  Así, las observaciones  pueden ser consideradas como funciones en los puntos de la geometría proyectiva.
\item
Si tenemos un objeto colocado en un plano distante y un lienzo en el cual vamos a dibujarlo, un dibujo en perspectiva se hace por la intersección de los puntos de la geometría proyectiva entre puntos del objeto real y puntos sobre el lienzo. Es decir un dibujo en perspectiva es simplemente una  homografía.
\item
El único invariante de la geometría proyectiva es la razón doble entre cuatro puntos distintos colineales y éste define todas  las reglas de la perspectiva.
\end{itemize}


\end{document}